\documentclass{amsart}
\pdfoutput=1
\usepackage{caption}
\usepackage{subcaption}
\usepackage{amssymb}
\usepackage{graphicx}
\usepackage{wrapfig}
\begin{document}

\author{G\'abor Fejes T\'oth}
\address{Alfr\'ed R\'enyi Institute of Mathematics,
Re\'altanoda u. 13-15., H-1053, Budapest, Hungary}
\email{gfejes@renyi.hu}

\title{Multiple arrangements}
\thanks{The English translation of the book ``Lagerungen in der Ebene,
auf der Kugel und im Raum" by L\'aszl\'o Fejes T\'oth will be
published by Springer in the book series Grundlehren der
mathematischen Wissenschaften under the title
``Lagerungen---Arrangements in the Plane, on the Sphere and
in Space". Besides detailed notes to the original text the
English edition contains eight self-contained new chapters
surveying topics related to the subject of the book but not
contained in it. This is a preprint of one of the new chapters.}

\begin{abstract}
This paper surveys the theory of multiple packings and coverings. The study
of multiple arrangements started in the 60s of the last century, and it
was restricted mostly to lattice arrangements on the plane or of general
arrangements of balls. We emphasize two new topics which were intensively
investigated recently: decomposition of multiple coverings into simple
coverings and characterization of multiple tilings.
\end{abstract}

\maketitle

We say that a family of sets is a {\it $k$-fold packing} if every point of the space
belongs to the interior of at most $k$ sets. Quite analogously, we say that a family
of sets forms a {\it $k$-fold covering} if every point of the space belongs to
the closure of at least $k$ sets. Let $\delta^k(K)$ and $\vartheta^k(K)$, respectively,
denote the densities of the densest $k$-fold packing and the thinnest $k$-fold
covering of the space with congruent copies of the convex body $K$. Similarly we use
the notation $\delta^k_T(K)$, $\vartheta^k_T(K)$, $\delta^k_L(K)$ and
$\vartheta^k_L(K)$ for the optimum densities of the corresponding $k$-fold translative
and lattice arrangements (compare the corresponding definitions of $\delta$, $\delta_T$,
$\delta_L$, $\vartheta$, $\vartheta_T$ and $\vartheta_L$ above).

\section{Multiple arrangements on the plane}

The literature on multiple packing and covering is relatively extensive,
and it mostly deals with arrangements of congruent copies of the circular disk
$B^2$. The values of $\delta^k(B^2)$ and
$\vartheta^k(B^2)$ are not known for any $k>1$. We know that
$\delta^2(B^2)>2\delta^1(B^2)$ and $\vartheta^2(B^2)<2\vartheta^1(B^2)$ as was shown by
{\sc{Heppes}}~\cite{Heppes55} and {\sc{Danzer}}~\cite{Danzer60}, respectively.
{\sc{G.~Fejes T\'oth}}~\cite{FTG76a} established the bounds
$$
\delta^k(B^2)\le k\,\frac{\pi}{6}\cot\frac{\pi}{6k}\quad{\rm{and}}\quad
\vartheta^k(B^2)\ge k\,\frac{\pi}{3}\csc\frac{\pi}{3k}.
$$
Observe that $\frac{\pi}{6}\cot\frac{\pi}{6k}$ and $\frac{\pi}{3}\csc\frac{\pi}{3k}$
are equal to the density of a disk with respect to the circumscribed and inscribed regular $6k$-gon.
For $k=1$ these inequalities are sharp as they coincide with Thue's and Kershner's theorems.

{\sc{Bolle}}~\cite{Bolle76} proved that there are positive constants $c_i$ such
that
$$k-c_1k^\frac{2}{5}\le\delta_L^k(B^2)\le{k}-c_2k^\frac{1}{4}$$
and
$$k+c_3k^\frac{1}{4}\le\vartheta_L^k(B^2)\le{k}+c_4k^\frac{2}{5}$$
and showed in~\cite{Bolle84} that the exponent $\frac{1}{4}$ is best possible in
these inequalities. In~\cite{Bolle89} he proved that for convex disks $K$
with piecewise twice differentiable boundary there are positive constants
$c(K)$ and $C(K)$ such that
$$\delta_L^k(K)\ge{k}-c(K)k^\frac{2}{5}$$
and
$$\vartheta_L^k(K)\le{k}+C(K)k^\frac{2}{5}.$$
Moreover, for a polygon $P$ the stronger inequalities
$$\delta_L^k(P)\ge{k}-c(P)k^\frac{1}{3}$$
and
$$\vartheta_L^k(P)\le{k}+C(P)k^\frac{1}{r3}$$
hold.

The exact values of $\delta^k_L(B^2)$ have been found for $k\le10$
(see {\sc{Heppes}}~\cite{Heppes59}, {\sc{Blundon}}~\cite{Blundon63},
{\sc{Bolle}}~\cite{Bolle76}, {\sc{Yakovlev}}~\cite{Yakovlev83}, {\sc{Temesv\'ari}}
~\cite{Temesvari94b} and {\sc{Te\-mes\-v\'ari}} and {\sc{V\'egh}}~\cite{TemesvariVegh}). The values of
$\vartheta^k_L(B^2)$ are known for $k\le8$ (see {\sc{Blundon}}~\cite{Blundon57}, {\sc{Haas}}~\cite{Haas},
{\sc{Subak}}~\cite{Subak} and {\sc{Temesv\'ari}}~\cite{Temesvari84a,Temesvari92a,Temesvari92b,Temesvari92a,Temesvari94a}.
{\sc{Linhart}}~\cite{Linhart83} described an algorithmic
approach for approximating the values of $\delta^k_L(B^2)$ with arbitrarily high
accuracy. Elaborating on results by {\sc{Yakovlev}}~\cite{Yakovlev84} {\sc{Temesv\'ari,
Horv\'ath}} and {\sc{Yakovlev}}~\cite{TemesvariHorvathYakovlev} described a method for finding the
densest $k$-fold lattice packing with circles. They reduced this task to a
finite number of optimization problems, each over an explicitly given compact
domain. A similar method for the thinnest $k$-fold lattice covering with circles
was given by {\sc{Temesv\'ari}}~\cite{Temesvari88}.

For a triangle $T$, {\sc{Sriamorn}}~\cite{Sriamorn14} determined $\delta_L^k(T)$ and
$\vartheta_L^k(T)$ for all $k$. Moreover, {\sc{Sriamorn}}~\cite{Sriamorn16}
showed that
$$\delta_T^k(T)=\delta_L^k(T)=\frac{2k^2}{2k+1}$$
and {\sc{Sriamorn}} and {\sc{Wetayawanich}}~\cite{SriamornWetayawanich15} showed that
$$\vartheta_T^k(T)=\vartheta_L^k(T)=\frac{2k+1}{2}$$
for all $k$.

It is worth mentioning that $\delta^k_L(B^2)=k\delta^1_L(B^2)$ for $k=2,\ 3$ and
$4$, and also $\vartheta^2_L(B^2)=2\vartheta^1_L(B^2)$. These equalities for the
very same multiplicities have been extended to an arbitrary centrally symmetric
convex disk in place of the circle by {\sc{Dumir}} and {\sc{Hans-Gill}}
\cite{DumirHans-Gill72a,DumirHans-Gill72b} and {\sc{G.~Fejes T\'oth}} \cite{FTG84a}.
The equality $\vartheta^2_L(B^2)=2\vartheta^1_L(B^2)$ was further generalized by
{\sc{Temesv\'ari}}, who proved in \cite{Temesvari84b} that the density of a 2-periodic
double covering by circles is at most $2\vartheta^1_L(B^2)$, and in \cite{Temesvari94c}
proved the analogous result for centrally symmetric convex disks. Recall from XX that
an ${m}$-periodic arrangement is the union of $m$ translates of a lattice arrangement.

\section{Decomposition of multiple arrangements}

The equalities $\delta^3_L(K)=3\delta^1_L(K)$ and $\delta^4_L(K)=4\delta^1_L(K)$
for centrally symmetric disks $K$ were derived by noticing that every $3$-fold lattice
packing by such a disk is the union of three simple ($1$-fold) packings, and every such
$4$-fold packing is the union of two $2$-fold packings. This observation belongs to
the topic concerning decompositions of multiple arrangements into simple ones,
problems and results that focus on the combinatorial structure of such arrangements.
Research in this direction was initiated by {\sc{Pach}}~\cite{Pach85}. He proved,
among other things, that every $2$-fold packing with positively homothetic copies of
a convex disk can be decomposed into four (simple) packings. For coverings, he made
the conjecture that for every convex disk $K$ there exists a minimal natural
number $m(K)$ such that every $m(K)$-fold covering of the plane by translates
of $K$ can be decomposed into two coverings. In \cite{Pach86} he proved this
conjecture for centrally symmetric polygons. New interest arose in the topic
after {\sc{Tardos}} and {\sc{T\'oth}}~\cite{TardosToth} proved the conjecture
for triangles. Soon after, {\sc{P\'alv\"olgyi}} and {\sc{T\'oth}}
~\cite{PalvolgyiToth} proved the conjecture for every convex polygon $P$.
Unfortunately, the number $m(P)$ increases with the number of sides of $P$,
thus the attempt to extend the result to all convex disks through polygonal
approximation fails. Still, it came as a surprise when {\sc{P\'alv\"olgyi}}
\cite{Palvolgyi} (see also {\sc{Pach}} and {\sc{P\'alv\"olgyi}}
\cite{PachPalvolgyi}) disproved Pach's conjecture by showing that it does
not hold for the circle. For subsequent developments on decomposition of
multiple arrangements we refer the reader to the survey article of
{\sc{Pach, P\'alv\"olgyi}} and {\sc{T\'oth}} \cite{PachPalvolgyiToth}.

\section{Multiple arrangements in space}

The densest 2-fold lattice packing and the thinnest 2-fold lattice covering of
balls in three dimensions were determined by {\sc{Few}} and {\sc{Kanagasabapathy}}
\cite{FewKanagasabapathy} and {\sc{Few}} \cite{Few67},
respectively. {\sc{Purdy}} \cite{Purdy} constructed a threefold lattice packing
of balls which he conjectured to be of maximum density. He supported the
conjecture by proving that it provides a local maximum of the density among
threefold lattice packings of balls.

Adapting Blichfeldt's idea, {\sc{Few}} \cite{Few64} gave the
following upper bound for the $k$-fold packing density of the $n$-dimensional
ball:
$$
\delta^k(B^n)\le (1+n^{-1})[(n+1)^k  -1][k/(k+1)]^{n/2}.
$$
This is better than the trivial bound $k$ only for large values of $n$ compared
to $k$. By a further elaboration on the same idea for $k=2$, {\sc Few}
\cite{Few68} obtained the stronger inequality
$$
\delta^2(B^n)\le\frac{4}{3}(n+2)\left(\frac{2}{3}\right)^{n/2}.
$$

{\sc G.~Fejes T\'{o}th} \cite{FTG79} gave a non-trivial upper bound for
$\delta^k(B^n)$, as well as a non-trivial lower bound for $\vartheta^k(B^n)$
for every $n$ and $k$.

For multiple lattice arrangements of balls, {\sc{Bolle}} \cite{Bolle79,Bolle82}
established sharper estimates. He proved that there are positive constants
$c_n$ and $C_n$ such that
$$
\frac{\delta_L^k(B^n)}{k}\le1-c_nk^{n+\frac{1}{n}}\quad{\hbox{and}}\quad
\frac{\vartheta_L^k(B^n)}{k}\ge1+c_nk^\frac{n+1}{2n}
$$
when $n\not\equiv1$ (mod 4)
and
$$
\frac{\delta_L^k(B^n)}{k}\le1-c_nk^\frac{n+3}{2n}\quad{\hbox{and}}\quad
\frac{\vartheta_L^k(B^n)}{k}\ge1+c_nk^\frac{n+3}{2n}
$$
when $n\equiv1$ (mod 4).

Extending the result of {\sc{Schmidt}} \cite{Schmidt61} to multiple
arrangements, {\sc{Florian}} \cite{Florian78a} proved that
$\delta^k(K)<k$ and $\vartheta^k(K)>k$ for every smooth
convex body $K$ without establishing a concrete bound.

By a Blichfeldt-type argument {\sc{Few}} \cite{Few64} proved
$$
\delta^k(B^n)\ge\delta(B^n)\left(\frac{2k}{k+1}\right)^{n/2}.
$$
In \cite{Few71} {\sc{Few}} studied the multiplicity of partial
coverings of space, and, as an application of a general theorem,
obtained a better lower bound for $\delta^k(B^n)$ for large values
of $k$ and $n$.

{\sc{Groemer}} \cite{Groemer86a} proved lower bounds for the $k$-fold
lattice packing density of a convex body $K$ involving the intrinsic
volumes of $K$. From his results follows the existence of positive
constants $c_n$ such that
$$\delta_L^k(K)\ge k-c_nk^{(n-1)/n}$$
for every convex body $K\in{E^n}$.

{\sc{Cohn}} \cite{CohnMJ} proved that
$$\vartheta_L^k(K)<[{(k+1)}^{1/n}+8n]^n=k(1+O(n^2k^{-1/n}))\qquad{\rm{as}}\ k\to\infty$$
for every $n$-dimensional convex body $K$.

The best known upper bound for $\vartheta^k(K)$ is due to {\sc{Nasz\'odi}}
and {\sc{Polyanskii}} \cite{NaszodiPolyanskii} who, improving slightly
on an earlier result by {\sc{Frankl, Nagy}} and {\sc{Nasz\'odi}}
\cite{FranklNagyNaszodi}, proved that
$$\vartheta^k(K)\le3.153(1+o(1))\max\{n\ln{n},k\}$$
for every convex body $K\in E^n$.

{\sc{Blachman}} and {\sc{Few}} \cite{BlachmanFew} gave bounds for the density
of multiple packings of spherical caps. Also {\sc L.~Fejes T\'oth} \cite{FTL66a},
{\sc{Galiev}} \cite{Galiev}, {\sc{Blinovsky}} \cite{Blinovsky} and {\sc{Blinovsky}}
and {\sc{Litsyn}} \cite{BlinovskyLitsyn} investigated multiple ball packings
in spherical spaces.

\section{Multiple tiling}

A system of bodies forms a {\it{$k$-fold tiling}} if each point of the space is
covered exactly $k$ times, except perhaps the boundary points of the bodies.
There are centrally symmetric polytopes that admit a translational
$k$-fold tiling, but no simple tiling. The simplest example is perhaps the
regular octagon of side-length 1, whose translates by the unit square
lattice form a 7-fold tiling.

{\sc{Bolle}} \cite{Bolle94} proved that a convex polygon that admits a
$k$-fold tiling of the plane by translations is centrally symmetric.
He also gave a characterization of those convex polygons that admit a
$k$-fold lattice-tiling. {\sc{Kolountzakis}} \cite{Kolountzakis21} gave an algorithm
which decides for a centrally symmetric convex polygon if it can tile the
plane by translations at some level. His algorithm runs in polynomial time
in the number of sides of the polygon. {\sc{Yang}} and {\sc{Zong}}
\cite{YangZong19} \cite{YangZong21} characterized those convex polygons that
admit two-, three-, four- or five-fold translational tiling. Only parallelograms
and centrally symmetric hexagons admit a two-, three- or four-fold translational
tiling. There are two more classes of polygons admitting five-fold tilings: the
affine images of a special octagon and of a decagon.

{\sc{Gravin, Robins}} and {\sc{Shiryaev}} \cite{GravinRobinsShiryaev}
proved that if translates of
a convex polytope form a $k$-fold tiling of $E^n$, then it is centrally
symmetric and its facets are centrally symmetric as well. This generalizes
a theorem of {\sc{Minkowski}} \cite{Minkowski97} concerning simple tilings.
For the three-dimensional case this means that only zonotopes admit a
translational $k$-fold tiling. For rational polytopes Gravin, Robins and
Shiryaev also proved the converse of their above mentioned theorem: Every
rational polytope in $E^n$ that is centrally symmetric and has centrally symmetric
facets admits a $k$-fold lattice tiling for some positive integer $k$.
For zonotopes, this was proved earlier by {\sc{Groemer}} \cite{Groemer78}.

A {\it{quasi-periodic}} set is a finite union of translated lattices, not
necessarily of the same lattice. {\sc{Kolountzakis}} \cite{Kolountzakis}
proved that a $k$-fold tiling by translates of a convex polygon other
than a parallelogram is quasi-periodic. {\sc{Gravin, Kolountzakis, Robins}}
and {\sc{Shiryaev}} \cite{GravinKolountzakisRobinsShiryaev} proved an
analogous theorem for the three-dimensional case: A $k$-fold tiling by
translates of a polytope that is not a two-flat zonotope is
quasi-periodic. A {\it{two-flat zonotope}} is the Minkowski sum of
two 2-dimensional symmetric polygons one of which may degenerate into
a single line segment.

{\sc{Gravin, Robins}} and {\sc{Shiryaev}} \cite{GravinRobinsShiryaev}
raised the problem whether the following generalization of the
Venkov-McMullen theorem holds: If translates of a polytope $P$ form a
$k$-fold tiling of $E^n$, then $P$ also admits an $m$-fold lattice
tiling for some, possibly different, multiplicity.  The
two-dimensional case of this conjecture was confirmed independently
by {\sc{Liu}} \cite{Liu} and {\sc{Yang}} \cite{Yang}. A further step in
the direction of proving the conjecture was made by {\sc{Chan}} \cite{Chan},
who proved it for certain quasi-periodic $k$-fold tilings. {\sc{Lev}}
and {\sc{Liu}} \cite{LevLiu} gave a characterization of those polytopes
in $E^n$ that tile with some multiplicity $k$ by translations along a
given lattice.

\small{
\bibliography{pack}}}
\mathrm{}